\documentclass[reqno]{amsart}
\usepackage{amsmath,amsfonts,amssymb}
\usepackage{verbatim}
\usepackage{enumerate}
\usepackage{url}
\def\N{\mathbb N}

\def\Q {\mathbb Q}

\def\F{\mathbb F}
\def\img{\mathop{\rm img}}

\def\rad{\mathop{\rm rad}}

\theoremstyle{plain}
\newtheorem{theorem}{Theorem}[section]
\newtheorem{lemma}[theorem]{Lemma}

\newtheorem{corollary}[theorem]{Corollary}

\newtheorem{remark}[theorem]{Remark}

\def\proof{{\it Proof: }}
\def\qed{\hfill\hbox{$\square$}}

\theoremstyle{definition}

\numberwithin{equation}{section}

\author[F.E. Brochero Mart\'{\i}nez]{F. E. Brochero Mart\'{\i}nez}

\address{
Departamento de Matem\'{a}tica\\
Universidade Federal de Minas Gerais\\
UFMG\\
Belo Horizonte, MG\\
 30123-970\\
 Brazil\\
 }
 \email{fbrocher@mat.ufmg.br }

\date{\today
}

\title{Structure of Finite dihedral 
group algebra}

\begin{document}
\maketitle

\begin{abstract} 
In this article, we show  the relation between the irreducible idempotents of  the cyclic group algebra $\F_qC_n$  and  the central irreducible  idempotents of the group algebras $\F_qD_{2n}$, 
where $\F_q$ is a finite field  with $q$ elements and $D_{2n}$ is the dihedral group of order $2n$, where $\gcd(q,n)=1$. 

In addition,  if every divisor of $n$ divides $q-1$, we  show  explicitly all central irreducible  idempotents of this group algebra and its Wedderburn decomposition.
\end{abstract}

\section{Introduction}
 Let $K$ be a field and $G$ be a group with $n$ elements.  It is known that, if   $char(K)\nmid n$, then the group algebra $KG$ is semisimple and as consequence of    Wedderburn Theorem, 
we have that  $KG$ is isomorphic to a direct sum of matrix algebras over division rings, such that each division algebra is a finite algebra over  the field $K$, i.e,  there exists an isomorphism 
$$\rho: KG\rightarrow M_{l_1}(D_1)\oplus  M_{l_2}(D_2)\oplus\cdots\oplus M_{l_t} (D_t),$$
where  $D_j$ are division rings such that $|G|=\sum_{j=1}^t l_j^2[D_j: K] 
$.  Observe that $KG$ has $t$  central irreducible idempotents, each one of the form 
$$e_i=\rho^{-1} (0\oplus\cdots\oplus 0\oplus I_i\oplus0\cdots\oplus 0),$$
where $I_i$ is the identity matrix of the component $M_{l_i}(D_j)$.  Then, the isomorphism $\rho$ determines explicitly each central irreducible idempotent. 

In the case $K=\Q$, the calculus of  central idempotents and Wedderburn decomposition is widely studied;  the classical method to calculate  the primitive central idempotents of group algebras depends on   computing   the  character group table.  
Other method is shown in \cite{JLP}, where  Jespers, Leal and Paques describe the central irreducible idempotents when $G$ is a nilpotent group, using the structure of its subgroups, without 
employing the characters of the group. 
Generalizations and improvements of this method can be found in  \cite{ORS}, where the authors   provide information about the Wedderburn decomposition of $\Q G$. 
This computational method is also used in \cite{BrRi} to compute the Wedderburn decomposition and the primitive
central idempotents of a semisimple finite group algebra  $KG$, where $G$ is  an abelian-by-supersolvable group
$G$  and  $K$ is a finite field. 

The structure of $KG$ when $G=D_{2n}$ is the dihedral group with $2n$ elements is well known for $K=\Q$ (see \cite{Gc}). 
In \cite{DFP}, Dutra, Ferraz and Polcino Milies impose conditions  over $q$ and $n$ in order for $\F_qD_{2n}$ to have the same  number of irreducible components that $\Q D_{2n}$.  This result is generalized in \cite{FGPM}, where Ferraz Goodaire and Polcino Milies find, for some families of groups,  conditions under $q$ and $G$ in order for $\F_qG$ to have the minimum number of simple components.

In this article, assuming that every prime factor of $n$ divides $q-1$, 
 we show explicitly  the central irreducible idempotents  of $\F_qD_{2n}$ and an isomorphism between the group algebra $\F_qD_{2n}$ and  its Wedderburn decomposition.
Observe that this isomorphism  also show the structure of $\mathcal U(\F_qD_{2n}) $, the unit group of $\F_qD_{2n}$.

\section{Idempotents of Cyclic Group Algebra}
Throughout this article, $\F_q$ denotes a finite field of order $q$, where $q$ is a power of a 
 prime and
$n$ is a positive integer such that
$\mbox{gcd}(n,q)=1$.
For every polynomial $g(x)$ with $g(0)\ne 0$,   $g^*$ denotes  the {\em reciprocal polynomial} of $g$, i.e., $g^*(x)=x^{\deg(g)}g(\frac 1x)$. The polynomial $x^n-1\in \F_q[x]$ splits in monic irreducible factors   as 
$$x^n-1=f_1f_2\cdots f_r f_{r+1}f_{r+1}^*f_{r+2}f_{r+2}^*\cdots f_{r+s}f_{r+s}^*,$$
 where $f_1=x-1$, $f_2=x+1$ if $n$ is even, and $f_j^*= f_j$  for $2\le j\le r$,  where $r$ is the number of auto-reciprocal factors in the factorization and $2s$ the number of non-auto-reciprocal factors.

We denote by  $C_n$ the cyclic group of order $n$. It is well known that $\F_qC_n\simeq \mathcal R_n=\frac {\F_q[x]}{\langle x^n-1\rangle }$, and by the Chinese  Remainder Theorem 
$$\frac {\F_q[x]}{\langle x^n-1\rangle}\simeq \bigoplus_{j=1}^{r+s} \frac  {\F_q[x]}{\langle f_j\rangle}\oplus \bigoplus_{j=r+1}^{r+s} \frac  {\F_q[x]}{\langle f^*_j\rangle}$$
is exactly the Wedderburn decomposition of the group algebra $\mathcal R_n$, so every primitive idempotent generates a maximal ideal of $\mathcal R_n$  and also one component of this direct sum.

In addition, since $\mathcal R_n$ is a principal ideal domain, every ideal of $\mathcal R_n$ is generated by a polynomial $g$ that is a divisor of $x^n-1$. The relation between the generator of the ideal and its principal idempotent is shown in the following lemma.

\begin{lemma}\label{idempotent1} Let $\mathcal I\subset \mathcal R_n$ be an ideal generated by the monic  polynomial $g$, that is divisor of $x^n-1$, and define   $f=\frac {x^n-1}{g}$.  Then the principal idempotent of $\mathcal I$ is
$$e_f=-\frac {((f^*)')^*}n  \cdot \frac {x^n-1}f
.$$
\end{lemma} 

\proof Let $t$ be an integer such that $n$ divides $q^t-1$. 
 By Theorem 2.1 in \cite{ArPr} (see also Theorem 3.4 in \cite{BrGi}), every primitive idempotent of $\frac{\F_t[x]}{\langle x^n-1\rangle}$ is given by 
$$u_\lambda=\frac {\lambda}n\cdot \frac{x^n-1}{x-\lambda}=\frac 1n\sum_{l=0}^{n-1} \lambda^{-l} x^l$$ 
where $\lambda^n=1$.

Since $f$ divides $x^n-1$, then $f$ splits in $\F_{q^t}[x]$ as $(x-\lambda_1)\cdots (x-\lambda_k)$ and
$$(f^*)'=\sum_{i=1}^k (-\lambda_i)\prod_{i\ne j} (1-\lambda_jx)=f^*\sum_{i=1}^k \frac {-\lambda_i}{x-\lambda_jx},$$
hence
$$e_f=-\frac {((f^*)')^*}n  \cdot \frac {x^n-1}f= \sum_{i=1}^k \frac{\lambda_i}{n}\cdot \frac {x^n-1}{x-\lambda}=\sum_{i=1}^k u_{\lambda_i}.
$$
Therefore $e_f$ is an idempotent of $\F_q[x]$.   In order to prove that $e_f$ is the principal idempotent of $\mathcal I$, it is enough to show that $g\cdot e_f=g$.  Observe that,  using partial fraction decomposition we obtain
$$g=\frac {x^n-1}{f}= \sum_{i=1}^k A_i u_{\lambda_i},$$
where $A_i=\frac n{\lambda_i} \frac 1{\prod_{j\ne i} (\lambda_i-\lambda_j)}$ and then
$$g\cdot e_f= \sum_{i=1}^k A_i u_{\lambda_i}\cdot\sum_{j=1}^k u_{\lambda_j}=\sum_{1\le i,j\le k} A_iu_{\lambda_i}u_{\lambda_j}= \sum_{i=1}^k A_i u_{\lambda_i}=g$$
as we wanted to prove.
\qed

\begin{remark} This lemma is also true for fields  with characteristic zero, it suffices to change in the proof the field $\F_{q^t}$ by the splitting field of the polynomial $f$.
\end{remark} 

\begin{corollary} The cyclic group ring $\mathcal R_n$ has $r+2s$ irreducible idempotents of the form $e_f$ given by Lemma \ref{idempotent1}, where  the polynomials $f$'s are  the  irreducible factors of $x^n-1\in \F_q[x]$.
\end{corollary}

\section{Central idempotents of Dihedral Group Algebra}

Throughout this section,   $\alpha_j$ denotes a root of the polynomial  $f_j$ and  
$D_{2n}$ denotes the dihedral group of order $2n$, i.e.
$$D_{2n}=\langle x,y| x^n=1,\, y^2=1, xy=yx^{-1}\rangle.$$
  We define integer numbers  $\epsilon$ and $\delta$ as
 $$\epsilon=\begin{cases}0 &\text{if char$(q)=2$}\\
 1&\text{if char$(q)\ne 2$ and $n$ is odd }\\
2&\text{if char$(q)\ne 2$ and $n$ is even}
\end{cases}$$ and $\delta=\max\{\epsilon, 1\}$.

The following theorem shows explicitly  the dependence of  the Wedderburn decomposition of the Dihedral group algebra over a finite field $\F_q$ with the factorization of $x^n-1\in \F_q[x]$.

\begin{theorem}\label{isomorphism}  The group algebra $\F_q  D_{2n}$ has  Wedderburn 
  decomposition  of the form
$$\F_q  D_{2n} \cong \bigoplus\limits_{j=1}^{r+s} A_j$$
where 
$$A_j= \begin{cases} \F_q \oplus \F_q 
&\text{if  $j\le \delta$,}\\
M_2(\F_q[\alpha_j+\alpha_j^{-1}])&\text{if $\delta+1\le j\le r$},\\
M_2(\F_q[\alpha_j])&\text{if $r+1\le j\le r+s$},\\
 
\end{cases} $$
\end{theorem}

\proof For each $j\in \{1,\dots, s+r\}$, let  $\tau_j$ be  the homomorphism  of $\F_q$-algebras defined by the generators of the group $D_{2n}$ as
$$\begin{array}{rccl}
\tau_1:&\F_qD_{2n}&\rightarrow&\F_q\oplus \F_q\\
&x&\mapsto&(1,1)\\
&y&\mapsto&(1,-1),
\end{array}
$$
in the case $\epsilon\ge 1$ and
$$\begin{array}{rccl}
\tau_2:&\F_qD_{2n}&\rightarrow&\F_q\oplus \F_q\\
&x&\mapsto&(-1,-1)\\
&y&\mapsto&(1,-1),
\end{array}
$$
in the case $\epsilon=2$, where the sum and product  in $\F_q\oplus \F_q$  is defined by adding and  multiplying   the corresponding components
of the same coordinates. Finally, for every $j\ge \epsilon+1$  
$$\begin{array}{rccl}
\tau_j:&\F_qD_{2n}&\rightarrow&M_2(\F_q[\alpha_j])\\
&x&\mapsto&\begin{pmatrix} \alpha_j &0\\ 0&\alpha_j^{-1}\end{pmatrix}\\
&y&\mapsto&\begin{pmatrix} 0 &1\\ 1&0\end{pmatrix}.
\end{array}
$$
It is easy to prove that $(\tau_j(x))^n=I$, $(\tau_j(y))^2=I$  and $\tau_j(x)\tau_j(y)=\tau_j(y)\tau_j(x)^{-1}$.

Observe that  in the case  of characteristic   $2$, i.e.  $\epsilon=0$, we have that 
$$\img(\tau_1)=\left\{\left. \begin{pmatrix} a&b\\ b&a\end{pmatrix}\right| a,b\in \F_q\right\},$$
that is isomorphic to $\F_q\oplus \F_q$ by the projection $\sigma_0$ over the first row of the matrix, 
where 
the product is defined by $(a,b)\cdot (c,d)=(ac+bd,ad+bc)$. Thus,
$\dim_{\F_q} (\img(\tau_1))=2$  in all cases. In addition, if  $n$ is even, then $\dim_{\F_q} (\img(\tau_2))=2$. 

For each $\delta<j\le r$, if  we define $Z_j=\begin{pmatrix} 1 &-\alpha_j\\ 1& -\alpha_j^{-1}\end{pmatrix}$, then 
$$\begin{array}{rccl}\sigma_j:&M_2(\F_q[\alpha_j]&\rightarrow &M_2(\F_q[\alpha_j]\\
&X&\mapsto&Z_j^{-1}XZ_j\end{array}$$
is an automorphism such that 
$$ \sigma_j\circ\tau_j(x)=\begin{pmatrix} 0&1\\ -1&\alpha_j+\alpha_j^{-1}\end{pmatrix}\quad\text{and}\quad \sigma_j\circ\tau_j(y)=\begin{pmatrix} 1&-(\alpha_j+\alpha_j^{-1})\\ 0&-1\end{pmatrix},
$$
so the images of the generators of $D_n$ are in $\F_q(\alpha_j+\alpha_j^{-1})$.
It follows that for each $j$ such that  $\delta< j\le r$ we have
$$\dim_{\F_q} (\img(\tau_j))=\dim_{\F_q} (\img(\sigma_j\circ\tau_j))\le 4\dim_{\F_q} (\F_q(\alpha_j+\alpha_j^{-1}))=2\deg(f_j)$$
and in the case  $r+1\le j\le r+s$, we know that
$$\dim_{\F_q} \img(\tau_j)\le 4\dim_{\F_q}\F_q(\alpha_j)=4\deg(f_j).$$

Now, let $\tau$ be the   homomophism  of $\F_q$-algebras defined by $\bigoplus\limits_{j=1}^{s+r} \tau_j$. Observe that this homomophism is injective. In fact, let $u$ be an element of $\F_qD_n$ such that $\tau(u)=0$. 
If we write  $u=P_1(x)+P_2(x)y$,  where $P_1$ and $P_2$ are polynomials of degree less than $n$, for each $j> \epsilon$, we have
$$\tau_j(u)=\begin{pmatrix}P_1(\alpha_j)&P_2(\alpha_j)\\ P_2(\alpha_j^{-1})&P_1(\alpha_j^{-1})\end{pmatrix}=\begin{pmatrix}0&0\\0&0\end{pmatrix},$$
so, $P_1(\alpha_j)=P_1(\alpha_j^{-1})=0$ and $P_2(\alpha_j)=P_2(\alpha_j^{-1})=0$.  In addition, if $\epsilon\ge 1$, then 
$$\tau_1(u)=(P_1(1)+P_2(1), P_1(1)-P_2(1)=0,$$
and if $\epsilon= 2$ we have
$$\tau_2(u)=(P_1(-1)+P_2(-1), P_1(-1)-P_2(-1)=0.$$
It follows that $P_1$ and $P_2$ are divisible by the polynomial $x^n-1$ and since the degrees of these polynomials are less that $n$, we conclude that  $P_1$ and $P_2$ are   null polynomials and therefore $\tau$ is an injective homomorphism. 

Finally, we observe that the homomorphism  $\rho:\F_qD_{2n}\to \bigoplus_{j=1}^{r+s} A_j$  
defined by $\rho=\bigoplus_{j=1}^{r+s} \rho_j$ where $\rho_j=\begin{cases} \sigma_j\circ \tau_j&\text{if $\epsilon<j\le r$}\\
\tau_j&\text{otherwise}\end{cases}$ is injective. Furthermore,   $\dim_{\F_q}(\F_qD_n)=2n$ and
\begin{align*}\dim_{\F_q}(\bigoplus_{j=1}^{r+s} A_j)&=2\epsilon+4\sum_{j=\epsilon+1}^{r} \dim_{\F_q}(\F_q[\alpha_j+\alpha_j^{-1}])+4\sum_{j=r+1}^{r+s} (\F_q[\alpha_j])\\
&=2\epsilon+2\sum_{j=\epsilon+1}^{r}\deg(f_j)+4\sum_{j=r+1}^{r+s} \deg(f_j)\\
&=2\deg(x^n-1)=2n.
\end{align*}
Therefore $\rho$ is an isomorphism.
\qed

\begin{remark} In the proof of the theorem we use the following facts:  if $\beta$ is a root of the polynomial $g\in \F_q[x]$, then $\beta^{-1}$ is root of the polynomial $g^*$. 
In addition,  when $g$ is auto-reciprocal and $\pm 1$ are not roots of $g$, there exists a polynomial $h\in \F_q[x]$ of degree $\frac{\deg(g)}2$, such that $\beta$ is a root of $g$ if and only if $\beta+\beta^{-1}$  is  a root of $h$.   
In fact, since $g$ is symmetrical, we can write $g$ as
$$g(x)=\sum_{j=0}^t a_j(x^{t+j}+x^{t-j})=x^t\sum_{j=0}^t a_j(x^{j}+x^{-j})=x^t\sum_{j=0}^t a_j D_j(z,1)=x^t h(z)$$
where $D_j$ is the Dickson polynomial of degree $j$ and $z=x+x^{-1}$ (see \cite{LMT} or \cite{Mey}).
\end{remark}

\begin{theorem}\label{idempotentdihedral} 
The dihedral group algebra  $\F_qD_{2n}$  has 
 $\epsilon+r+s$ central irreducible idempotents:
\begin{enumerate}
\item $2\epsilon$ idempotents of the form $\frac {1+y}2 e_{f_j}$ and $\frac {1-y}2 e_{f_j}$, where $j\le \epsilon$.
\item $r-\epsilon$ idempotents  $e_{f_j}$, where $j=\epsilon+1, \dots, r$,  generated by the auto-reciprocals factor of $x^n-1$. 
\item $s$ idempontents  $e_{f_j}+e_{f_j^*}$, where $j=r+1,\dots,r+ s$.
\end{enumerate}

\end{theorem}

\proof
Since the homomorphism $\tau$  in the proof of Theorem \ref{isomorphism}  is injective, then the image of a central primitive  idempotent  $u$ by the homomorphism has to be zero in every component, except for one   component  where the image is the identity, i.e., 
for some $i$ fixed, $\tau_j(u)=\delta_{i,j} I_j$, where $I_j$ is the identity over the component $A_j$. 
Let  $u=P(x)+Q(x)y$ be a representation of $u$, where $P$ and $Q$ are polynomials in $\F_q[x]$ of degree less than or equal to $n-1$.  Observe that $Q$ is zero when  calculated  at each root of the polynomial $x^n-1=0$, so  $Q$ is the null polynomial. 
In addition, $P$ is one when we calculate it at the roots of the polynomials $f_j$ and $f_j^*$ and zero when we calculate  it at the other roots of the polynomial $x^n-1$.
The unique polynomial of degree less or equal to $n-1$  that satisfies that proprieties is $e_{f_j}$, when  $f_j=f_j^*$  and 
 $e_{f_j}+e_{f_j^*}$ when $f_j\ne f_j^*$.
Finally,  if $j\le \epsilon$ the image $\tau_j(e_{f_j})=(1,1)$  is not a primitive idempotent, and we can decompose this idempotent in two central primitive idempotents, $(\frac {1+y}2)e_{f_j}$ and $(\frac {1-y}2)e_{f_j}$, such that $\tau_j((\frac {1+y}2)e_{f_j})=(1,0)$ and $\tau_j((\frac {1-y}2)e_{f_j})=(0,1)$.
\qed

\section{Explicit form of the Idempotents when $\rad(n)|(q-1)$}
Throughout this section,  we assume that every prime factor of $n$ divides $q-1$, $\kappa$ and $\nu$ denote the numbers $\gcd(n,q-1)$ and $\min\{\nu_2(\frac n2),\nu_2(q+1)\}$  respectively, $\theta$  and $\alpha$ are   generators of $\F_q^*$ and  $\F_{q^2}^*$  such that $\alpha^{q+1}=\theta$.
In the following results, we  show   the explicit form of the idempotents of the cyclic group algebra $\F_qC_n$ and   the Wedderburn decomposition of  the Dihedral group algebra  $\F_qD_{2n}$
.  In other to show that representation, we need the following lemma

\begin{lemma}\cite[Corollary 3.3 and Corollary 3.6]{BGO} \label{factorsb}
The factorization of $x^n-1$ in irreducible factors of $\F_q[x]$ depends on $n$ and $q$ in the following form:
\begin{enumerate}[(i)]
\item If $8\nmid n$ or $q\not\equiv 3 \pmod 4$, then
$$x^n-1=\prod_{t|m}\prod_{{1\le u\le \gcd(n,q-1)\atop \gcd(u,t)=1}} (x^t-\theta^{ul}),$$
where   $m=\frac n{\kappa}$ and $l=\frac {q-1}{\kappa}$. 
In addition,  for each  $t$ such that $t|m$, the number of irreducible factors of degree $t$ is $\frac{\varphi(t)}t\cdot \kappa$,  where $\varphi$ denotes the Euler Totient function.
\item If $8\mid n$ and $q\equiv 3 \pmod 4$, then 
$$x^n-1=\prod_{t|m'\atop t\text{ odd}}\prod_{{1\le w\le \kappa}\atop \gcd(w,t)=1} (x^t-\theta^{wl})\cdot \prod_{t|m'}\prod_{ u\in\mathcal S_t} (x^{2t}-(\alpha^{ul'}+\alpha^{qul'})x^t+\theta^{ul'}),$$
where 
$m'=\frac n{2^{\nu}\kappa}$,  $l'=\frac {q^2-1}{2^{\nu}\kappa}$, 
 and   $\mathcal S_t$ is the set 
 $$\left\{u\in\N\left|  {1\le u\le 2^{\nu}\kappa, \gcd(u,t)
=1\atop 
2^{\nu}\nmid u \text{ and }\ u< \{qu\}_{2^r\kappa}}\right.\right\},$$
where $\{a\}_b$ denotes the  remainder of the division of $a$ by $b$, i.e.   the number $0\le c<b$ such that $a\equiv c\pmod b$.
In addition, for each $t$  odd such that $t|m'$, the number of irreducible binomials  of degree $t$ and   $2t$ is
$\dfrac{\kappa\cdot\varphi(t)}t$ and
$\dfrac{\kappa\cdot\varphi(t)}{2t}$    respectively,  and
the number of irreducible trinomials  of degree $2t$ is
$$
\begin{cases}
\dfrac{\varphi(t)}t\cdot 2^{\nu-1} \kappa,&\text{if $t$ is even}\\
\dfrac{\varphi(t)}{t}\cdot (2^{\nu-1}-1) \kappa,&\text{if $t$ is odd}.
\end{cases}
$$
\end{enumerate}
\end{lemma}

The following corollary, direct from   Lemmas  \ref{idempotent1}  and \ref{factorsb}, shows the explicit form of each idempotent  of the cyclic group algebra $\F_qC_n$ when $\rad(n)|(q-1)$.

\begin{corollary} \label{idemR_n} Let $m$, $m'$, $l$ and $l'$ be as in  Lemma \ref{factorsb}.
\begin{enumerate}
\item If $8\nmid n$ or $n\not\equiv 3 \pmod 4$, then every irreducible  idempotent of the ring $\mathcal R_n$  is of the form
$$e_{t,ul}=\frac {\theta^{ul}t}n\cdot\frac{x^n-1}{x^t-\theta^{ul}}$$
where $t$ and $u$  satisfy  the condition of Lemma \ref{factorsb} item {\em (i)}.
\item If $8|n$ or $n\equiv 3 \pmod 4$, then every irreducible   idempotent of the ring $\mathcal R_n$  is of the form shown in  {\em (1)} and  of the form 
$$ e_{t, ul'}= \frac{t}{n} \left((\alpha^{ul'}+\alpha^{ul'q})x^t - 2\theta^{ul'} \right)\frac{x^n-1}{(x^{2t}-(\alpha^{ul'}+\alpha^{qul'})x^t+\theta^{ul'}}
,$$
where $t$ and $u$  satisfy  the condition of Lemma \ref{factorsb} item {\em (ii)}.
\end{enumerate}
\end{corollary}

\begin{remark} By Theorem \ref{idempotentdihedral} If 
\begin{itemize}
\item   char$(\F_q)=2$, or 
\item $n$ is odd and $\theta^{ul}\ne1$,  or 
\item $n$ is even and $\theta^{ul}\ne\pm 1$,  
\end{itemize}
then every idempotent found  in Corollary \ref{idemR_n} item (1) is   a central  irreducible idempotent of $\F_qD_{2n}$.  Otherwise, the idempotent can be reduced to two central primitive idempotents $\frac {1+y}2e_{t,ul}$ and $\frac {1-y}2e_{t,ul}$.

In addition, $e_{t,ul'}$ of item (2) is also a central irreducible idempotent of $\F_qD_{2n}$  if
$\theta^{ul'}= 1$,  otherwise, the central irreducible idempotent is $e_{t,ul'}+ e_{t,-ul'}$.
\end{remark}

\begin{theorem} 
The Wedderburn decomposition of  the group algebra  $\F_qD_{2n}$ depends on $n$ and $q$ in the following form:
\begin{enumerate}
\item When $n$ is odd, the decomposition is
$$2\F_q\oplus\frac {\kappa-1}2 M_{2}(\F_q)\oplus\bigoplus_{t|m \atop t\ne 1}  \frac {\kappa\cdot\varphi (t)}{2t} M_{2}( \F_{q^t}).$$
\item When $n$ is even,
\begin{enumerate}[(2.1)]
\item if $q\equiv 1 \pmod 4$ or $8\nmid n$, the decomposition is
$$4\F_q\oplus\left(\frac {\kappa}2-1\right) M_{2}(\F_q)\oplus\bigoplus_{t|m \atop t\ne 1}  \frac {\kappa\cdot\varphi (t)}{2t} M_{2}( \F_{q^t}),$$
\item  if $q\equiv 3 \pmod 4$ and $8| n$, the decomposition is
$$4\F_q\oplus(\kappa+2^{\nu-i}-3) M_{2}(\F_q)\oplus(
2^{\nu-2}\kappa-2^{\nu-1}-\frac k4+1
) M_{2}(\F_{q^2})\oplus\bigoplus_{{t|m' \atop t\text{ odd}}\atop t\ne 1}  \frac {\kappa\cdot\varphi (t)}{2t}  M_{2}( \F_{q^t})
$$
$$
\oplus\bigoplus_{t|m' \atop t\text{ even}} 
\dfrac{2^{\nu-2}\kappa\cdot\varphi(t)}t  M_{2}( \F_{q^{2t}})
\oplus\bigoplus_{{t|m' \atop t\text{ odd}}\atop t\ne 1} 
\dfrac{ (2^{\nu-1}-1) \kappa\cdot\varphi(t)}{2t} M_{2}( \F_{q^{2t}}).
$$
where $i=\begin{cases}0&\text{if $\nu_2(q+1)> \nu_2(\frac n2)$}\\1&\text{if $\nu_2(q+1)\le \nu_2(\frac n2)$.}\end{cases}$
\end{enumerate}
\end{enumerate}
\end{theorem}
\proof
First, we consider the case $n\not\equiv 3\pmod 4$ or $8\nmid n$, so  every irreducible factor of $x^n-1$ is a binomial,  and except for the factors $x-1$ and $x+1$,  we have that  any irreducible factor of the form $x^t-a$ is  not auto-reciprocal.  Thus, we have two cases to analyse:

\begin{enumerate}[$i)$]
\item
If $n$ is odd, we have  that $\epsilon=0$ or $1$  and $r=1$. By Lemma \ref{factorsb} there exist $\frac {\kappa\varphi(t)}t$ irreducible factors of degree $t$  and  by Theorem  \ref{idempotentdihedral}
 there exist two components isomorphic to $\F_q$, $\frac {\kappa\varphi(t)}{2t}$ components of the form $M_2(\F_{q^t})$ if $t>1$ and 
$\frac {\kappa-1}2$ components of the form $M_2(\F_q)$ if $t=1$, where $t$ is a divisor of $m$.  So we obtain  item (1).

\item If $n$ is  even, we have  that $\epsilon=2$ and there exist four components isomorphic to $\F_q$.  
In addition,   every factor of $x^n-1$  different that $x\pm 1$ is  a non-auto-reciprocal   binomial, then $r=2$, and by the same argument of the previous case there exist   $\frac {\kappa\varphi(t)}{2t}$ components of the form $M_2(\F_{q^t})$ if $t>1$ and 
$\frac {\kappa-2}2$ components of the form $M_2(\F_q)$ if $t=1$, where $t$ is a divisor of $m$. So, we obtain  item (2.1).
\end{enumerate}

Finally, in the case which  $q\equiv 3 \pmod 4$ and $8| n$, every factor of $x^n-1$  is a binomial or a trinomial.  The unique auto-reciprocal factor of the form $x^t-a$ with $t$ odd is $f_1=x-1$. 
Now, suppose that $x^{2t}-(\alpha^{ul'}+\alpha^{qul'})x^t+\theta^{ul'})$ is an irreducible factor of $x^n-1$ as in Lemma  \ref{factorsb} item (b), such that  it is an auto-reciprocal polynomial. It follows that $\theta^{ul'}=1$ and therefore $(q-1)|ul'$. Since 
$$l'=\frac{q-1}{\gcd(n,q-1)}\cdot \frac{q+1}{2^\nu},$$
the polynomial is auto-reciprocal when $\gcd(n,q-1)| u\cdot\frac{q+1}{2^\nu} $ and  we have two cases to consider:
\begin{enumerate}[$i)$]
\item If $\nu_2(q+1)\le \nu_2(\frac n2)$ then $\frac {q+1}{2^\nu}$ is odd and $\gcd (\gcd(n,q-1), \frac {q+1}{2^\nu})=1$, therefore $\gcd(n,q-1)|u$. But $t|m'|n$ and $\gcd(t,u)=1$, then these conditions imply that $t=1$ and $u$ is a multiple of $\gcd(n,q-1)$ not divisible by $2^\nu$ and less than $2^\nu\gcd(n,q-1)$. So there exist $2^{\nu} -2$  values of $u$ that generate $2^{\nu-1}-1$ auto-reciprocal factors, all of them of degree $2$, each one generating a component of the form $M_2(\F_q)$.
 In addition, we have $\kappa-2$ irreducible factors of degree $1$, each one generating a component of the same type.

Therefore there exist $(\kappa-2)+(2^{\nu-1}-1)=\kappa+2^{\nu-1}-3$ components $M_2(\F_q)$ and $$\frac k4(2^\nu-1)-(2^{\nu-1}-1)=2^{\nu-2}\kappa-2^{\nu-1}-\frac k4+1$$ components $M_2(\F_{q^2})$.

\item  If $\nu_2(q+1)> \nu_2(\frac n2)$ then $\frac {q+1}{2^\nu}$ is even and $\gcd (\gcd(n,q-1), \frac {q+1}{2^\nu})=2$, therefore $\frac 12\gcd(n,q-1)|u$. Similarly, we obtain $t=1$ and $u$ is a multiple of $\frac 12\gcd(n,q-1)$ non divisible by $2^\nu$ and less than $2^\nu\gcd(n,q-1)$. So there exist $2^{\nu+1} -2$  values of $u$ and then $2^\nu-1$  auto-reciprocal factors, all of them of degree $2$.

Then there exist $\kappa+2^{\nu}-3$ components $M_2(\F_q)$ and $2^{\nu-2}\kappa-2^{\nu}-\frac k4+1$
components $M_2(\F_{q^2})$.\qed
\end{enumerate}


\begin{thebibliography}{99}
\bibitem{ArPr} Arora, S. K., Pruthi M. {\em Minimal Codes of Prime-Power Length}. Finite Fields Appl. {\bf 3} (1997) 99-113

\bibitem{BrRi}  Broche, O.,  del R\'io, \'A., {\em Wedderburn decomposition of finite group algebras, Finite Fields Appl.} {\bf 13} (2007) 71-79.

\bibitem{BrGi}  Brochero Mart\'inez, F.E., Giraldo Vergara, C.R., {\em Explicit Idempotents of Finite Group Algebra} Finite Fields Appl. {\bf 28} (2014)  123-131

\bibitem{BGO} Brochero Mart\'inez, F.E., Giraldo Vergara, C.R., Batista de Oliveira, L., {\em Explicit Factorization of $x^n-1\in \F_q[x]$},  submitted for publication in Designs, Codes and Cryptography.  Preprint available on \url{http://arxiv.org/abs/1404.6281}


\bibitem{DFP} Dutra,F., Ferraz, R., Polcino Milies, C., {\em Semisimple group codes and dihedral codes},  Algebra Discrete Math. {\bf 3} (2009),  28-48.
\bibitem {FGPM} Ferraz, R., Goodaire, E., Polcino Milies, C., {\em Some classes of semisimple group (and loop) algebras over
finite fields}, J. of Algebra {\bf 324} (2010) 3457-3469

\bibitem{Gc} 
Giraldo Vergara, C.R.,  Brochero Mart\'inez, F. E.
Wedderburn decomposition of some special rational group algebras. 
Lect. Mat. 23 , no. 2, (2002)  99-106.





\bibitem{JLP}Jespers, E.,   Leal, G.,  Paques, A.,  {\em Central idempotents in rational group algebras of finite nilpotent groups}, J. Algebra Appl. 2 (1) (2003) 57-62.

\bibitem{LMT}  Lidl, R.  Mullen, G. L., Turnwald, G.,  {\em Dickson polynomials}, Pitman Monographs and
Surveys in Pure and Applied Math., Longman, London-Harlow-Essex, 1993.


\bibitem{Mey} Meyn H.{\it Factorization of the cyclotomic polynomials $x^{2^n}+1$ over finite fields}. Finite Fields Appl. 2, (1996) 439-442


\bibitem{ORS}  Olivieri, A., del R\'io, \'A,  Sim\'on, J.J., {\em On monomial characters and central idempotents of rational group algebras}, Comm. Algebra {\bf 32} (4) (2004) 1531-1550.


\end{thebibliography}
\end{document}